\documentclass{amsart}
\usepackage{amsmath,amsthm,amssymb}

\usepackage{eucal}
\usepackage{verbatim}
\usepackage[dvips]{graphicx}
\usepackage{multirow}
\usepackage{fancyhdr}
\usepackage{color}
\usepackage{enumerate}
\usepackage{calrsfs}
\usepackage{fullpage}

\numberwithin{equation}{section}

\newtheorem{theorem}{Theorem}
\newtheorem{lemma}{Lemma}

\newtheorem{proposition}{Proposition}

\newtheorem*{convention}{Convention}

\makeatletter
\newcommand{\leqnomode}{\tagsleft@true}
\newcommand{\reqnomode}{\tagsleft@false}
\makeatother

\def\({\begin{eqnarray}}
\def\){\end{eqnarray}}
\def\[{\begin{eqnarray*}}
\def\]{\end{eqnarray*}}
\def\part#1#2{\frac{\partial #1}{\partial #2}}

\def\R{\mathbb{R}}
\def\d{\mathrm{d}}
\def\tot#1#2{\frac{\d #1}{\d #2}}
\def\eps{\varepsilon}

\def\L{\mathcal{L}}
\def\E{\mathbb{E}}

\def\tw{\widetilde w}
\def\ty{\widetilde y}
\def\L{\mathcal{L}}

\begin{document}

\title{Asymptotic and exponential decay in mean square for delay geometric Brownian motion}   

\author{Jan Haskovec${}^1$}\footnote{Computer, Electrical and Mathematical Sciences \& Engineering, King Abdullah University of Science and Technology, 23955 Thuwal, KSA.
jan.haskovec@kaust.edu.sa}

\maketitle

\begin{abstract}
We derive sufficient conditions for asymptotic and monotone exponential decay
in mean square of solutions of the geometric Brownian motion with delay.
The conditions are written in terms of the parameters and are explicit
for the case of asymptotic decay.
For exponential decay, they are easily resolvable numerically.
The analytical method is based on construction of a Lyapunov functional
(asymptotic decay) and forward-backward estimate for the square mean (exponential decay).
\end{abstract}
\vspace{2mm}

\vspace{2mm}

\textbf{2010 MR Subject Classification}: 34K11, 34K25, 34K50, 60H10
\vspace{2mm}

\section{Introduction and main result}\label{sec:Intro}
Geometric Brownian motion (also called Ornstein-Uhlenbeck process with multiplicative noise)
is the strong solution of the It\^{o} stochastic differential equation
\( \label{eq:GBM}
   \d w(t) = - \lambda w(t) \d t + \sigma w(t) \d B_t,
\)
where $\lambda$ and $\sigma$ are real parameters
and $\d B_t$ denotes the one-dimensional Wiener process.
It is one of the stochastic processes very often used in
applications, in particular in financial mathematics
to model stock prices in the Black-Scholes model \cite{Hull}.
However, modelling the price process by geometric Brownian
motion has been criticized because the past of the volatility is not taken into account.
Consequently, \cite{Appleby} suggests to replace the multiplicative
constants $\lambda$ and $\sigma$ in \eqref{eq:GBM} by some linear functionals on the
space of continuous functions.
Here we make the generic choice of constant delay model, i.e.,
we evaluate $w$ in the right-hand side of \eqref{eq:GBM}
at the past time instant $t-\tau$, with $\tau>0$.
This leads to the following delay It\^{o} stochastic differential equation
\( \label{eq:0}
   \d w(t) = - \lambda w(t-\tau) \d t + \sigma w(t-\tau) \d B_t.
\)
The main goal of this paper is to derive sufficient conditions
for asymptotic and monotone (exponential) decay
in mean square of the solutions of \eqref{eq:0}.


Solutions of delay (retarded) differential equations are well known to develop
oscillations in certain regimes \cite{Smith}.
Taking the expectation of \eqref{eq:0}, we obtain the deterministic
delay differential equation for $u(t):=\E[w(t)]$,
\(  \label{eq:u}
    \dot u(t) = -\lambda u(t-\tau).
\)
Despite its simplicity, it exhibits a surprisingly rich qualitative dynamics.
An analysis of the corresponding characteristic equation
\[
   z + \lambda\tau e^{-z} = 0,
\]
where $z\in\mathbb{C}$, reveals that:
\begin{itemize}
\item
If $0 < \lambda\tau < e^{-1}$, then $u=0$ is asymptotically stable.
Solutions of \eqref{eq:u} subject to constant nonzero initial datum on $[-\tau,0]$
tend to zero monotonically (exponentially) as $t\to\infty$.
\item
If $e^{-1} < \lambda\tau < \pi/2$, then $u=0$ is asymptotically stable,
but every nontrivial solution of \eqref{eq:u} is oscillatory,
i.e., changes sign infinitely many times on $(0,\infty)$.
\item
If $\lambda\tau > \pi/2$, then $u=0$ is unstable.
\end{itemize}
We refer to Chapter 2 of \cite{Smith} and \cite{Gyori-Ladas} for details.
Consequently, two very natural questions arise in connection
with (linear) delay differential equations: Under which conditions
does the solution tend to zero asymptotically as $t\to\infty$,
and under which conditions is this decay monotone?
This paper is devoted to the study of these two questions
in mean square sense for solutions of \eqref{eq:0}.

Various types of sufficient conditions for stability (in some sense)
of equation \eqref{eq:0} and its generalizations have been established
in the literature, see \cite{KM99, Mao2, Mohammed} for an overview.
However, to our best knowledge, none of them provide
an explicit formula relating the parameters $\lambda$, $\tau$ and $\sigma$.
A remarkable result by~\cite{Appleby} states that
$$
\lim_{t\to\infty} 
\E \left[|w(t)|^2 \right] = 0
\qquad
\qquad 
\mbox{if and only if}
\qquad
\qquad 
\int_0^\infty r_\lambda(t)^2 \d t < \frac{1}{\sigma^2},
$$
where $w$ is a solution of \eqref{eq:0} and
$r_\lambda$ is the fundamental solution of the delayed ODE \eqref{eq:u},
i.e., formally, $r_\lambda$ solves~\eqref{eq:u} subject to the 
initial condition $u(t) = \chi_{\{0\}}(t)$ for $t\in(-\tau,0]$.
The fundamental solution $r_\lambda$ can be constructed by the 
method of steps~\cite{Smith}, however, to out best knowledge,
analytic evaluation of its  $L^2(0,\infty)$-norm is an open problem.
An explicit sufficient condition for asymptotic mean square stability
of \eqref{eq:0} has been provided in \cite{EHS}, together with
numerical experiments (systematic Monte Carlo simulations)
giving a hint about how far the analytical result is from optimal.
However, \cite{EHS} considers \eqref{eq:0} only as a special case
of a more general delay stochastic system, which leads to some
inefficiencies. Our first result, Theorem \ref{thm:asymptotic},
improves the sufficient condition of \cite{EHS},
and is still explicit in terms of the parameter values.
The proof is based on a construction of an appropriate Lyapunov functional.
Our second result, Theorem \ref{thm:exponential},
is based on a forward-backward estimate for the mean square and
provides sufficient condition for exponential decay
in mean square of solutions of \eqref{eq:0}.
The condition, written in terms of $\lambda$, $\tau$ and $\sigma$,
is not fully explicit, however, can be very easily resolved 
numerically.

This paper is organized as follows.
In Section \ref{sec:main} we provide an overview
of our results, formulate the corresponding theorems
and discuss their optimality. 
In Section \ref{sec:asymptotic} we provide the proof
for the case of asymptotic decay, which is based
on a construction of an appropriate Lyapunov functional.
In Section \ref{sec:exponential} we provide the
proof of exponential decay, based on forward-backward
estimates for the mean square of the solution.

\section{Main results}\label{sec:main}

A simple scaling analysis of \eqref{eq:0} reveals that
its dynamics depends on two parameters, which can be chosen
as $\lambda\tau$ and $\sigma/\sqrt{\lambda}$.
Therefore, with abuse of notation, we rename
$\lambda\tau \mapsto \tau$ and $\sigma/\sqrt{\lambda} \mapsto \sigma$
and rewrite \eqref{eq:0} as
\(  \label{eq:main}
    \d w(t) = - w(t-\tau) \d t + \sigma w(t-\tau) \d B_t.
\)
We shall consider \eqref{eq:main} subject to the deterministic initial datum
\(   \label{IC:main}
   w(s) = w_0(s) \qquad\mbox{for } s \in [-\tau,0],
\)
where $w_0=w_0(s)$ is a continuous function on $[-\tau,0]$.
We have the following result regarding the well posedness
of the problem \eqref{eq:main}--\eqref{IC:main}.

\begin{proposition}\label{thm:ex}
The stochastic delay differential equation~\eqref{eq:main} with initial
datum~\eqref{IC:main} admits a unique global solution 
$w=w(t)$ on $[-\tau,\infty)$ 
which is an adapted process with 
$\E\left[\int_{-\tau}^T |w(t)|^2 \d t\right] < \infty$ for all $T<\infty$.
\end{proposition}

\proof
The proof follows directly from Theorem 3.1 of \cite{Mao} and the subsequent remark on p. 157 there.
In particular, the right-hand side of \eqref{eq:main} is independent of the present state $w(t)$,
so that the solution can be constructed by the method of steps~\cite{Smith}.
The second order moment is bounded on any bounded interval
due to the linearity of the equation.
\endproof

\medskip

\begin{convention}
Throughout the paper we adopt the following notational convention:
we denote $\tw$ the quantity $w$ evaluated at time $t-\tau$, i.e., $\tw:=w(t-\tau)$,
while $w$ shall denote $w:=w(t)$. The same convention shall be applied to any other
time-dependent variable, in particular, the quantity $y:=\E[w^2/2]$ that we shall use
in the sequel.
\end{convention}
\medskip

Our first result gives an explicit sufficient condition in terms of
the parameters $\tau$, $\sigma$
for asymptotic decay of the square mean $\E[w(t)^2]$
for solutions $w=w(t)$ of \eqref{eq:main}.

\begin{theorem}\label{thm:asymptotic}
Let 
\(  \label{cond:asympt}
    \sigma^2<2,\qquad
    \tau < 1 - \sqrt{\sigma^2-\frac{\sigma^4}{4}}.
\)
or, equivalently,
\(  \label{cond:asympt2}
   \tau<1,\qquad
    \sigma < \sqrt{2-\tau}  -\sqrt{\tau}.
\)
Then the solutions $w=w(t)$ of \eqref{eq:main} satisfy
\[
   \lim_{t\to\infty} \E[w(t)^2] = 0.
\]
\end{theorem}

Let us observe that the above result is suboptimal
in the borderline case $\sigma=0$,
i.e., the deterministic regime given by \eqref{eq:u}.
Indeed, \eqref{cond:asympt} then turns into $\tau < 1$,
while solutions of \eqref{eq:u}
asymptotically decay to zero if (and only if) $\tau < \pi/2$, see, e.g., \cite{Smith}.
However, in the other borderline case $\tau=0$,
\eqref{cond:asympt2} becomes $\sigma^2<2$,
which is the sharp condition for asymptotic vanishing
of the mean square of geometric Brownian motion
\eqref{eq:GBM}, see, e.g., \cite{Oksendal, Mao}.

\begin{figure}
\centerline{
\includegraphics[width=0.6\columnwidth]{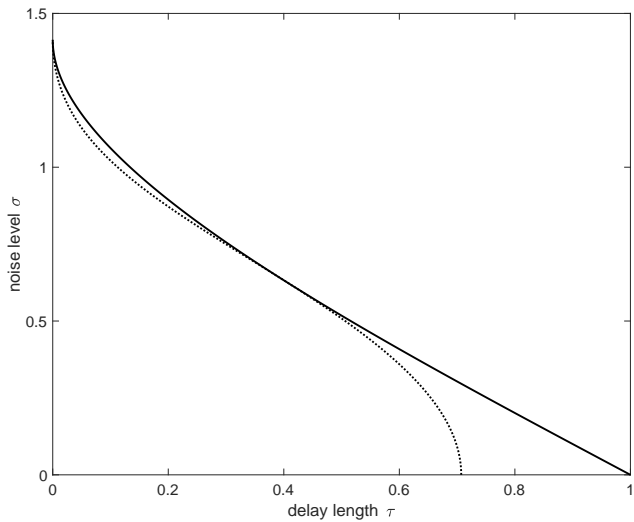}}
\caption{Comparison of sufficient conditions for asymptotic decay of $\E[w^2]$ 
formulated by Theorem \ref{thm:asymptotic}, solid line,
and by \cite[Lemma 3.5]{EHS}, dashed line.}
\label{fig:1}
\end{figure}

We also note that the result of \cite{EHS}
provides a less optimal condition than Theorem \ref{thm:asymptotic}.
Indeed, the condition stated by Lemma 3.5 of \cite{EHS} reads, in our notation,
\(   \label{cond:EHS}
   \sigma^2 < 2,\qquad
   \tau < \frac{1}{4} \left( -2\sigma^2 + \sqrt{4\sigma^4 + 2(2-\sigma^2)^2} \right).
\)
As illustrated in Fig. \ref{fig:1}, the upper bound on $\tau$ of \eqref{cond:EHS}
is more restrictive then the one of \eqref{cond:asympt} for all values of $\sigma^2<2$.
We see that Theorem \ref{thm:asymptotic} represents an improvement 
especially in the low noise regime. 
In the limit $\sigma^2\to 0$ it improves the restriction $\tau<1/\sqrt{2}$ imposed by \eqref{cond:EHS}
to $\tau<1$ (which, however, is still not optimal, as noted above).

Finally, let us refer to \cite[Fig. 2]{EHS} for a comparison of the analytical condition \eqref{cond:EHS}
to results of systematic Monte Carlo simulations, which indicates that there is still a significant
potential for improvement of the analytical result.

Our second result provides a sufficient condition
for exponential (monotone) decay of the square mean 
of solutions of \eqref{eq:main}.
Obviously, monotonicity of the solution strongly depends on the initial datum $w_0$.
Therefore, we consider the generic case of constant, nonzero initial condition $w_0\in\R$ in the below Theorem.
For notational convenience we define, for $s\geq 0$, the function $G=G(s)$,
\( \label{def:G}
   G(s) := \left( \frac{e^{2s} - 1}{2s} \right)^{1/2}.
\)
\begin{theorem}\label{thm:exponential}
Let $\tau$, $\sigma\geq 0$ be such that 
the conditions
\(  \label{ass:sigma1}
   \sigma < e^{-\mu\tau} \sqrt{2\mu - 2e^{\mu\tau}}
\)
and
\(  \label{ass:sigma2}
   \sigma < -G(\mu\tau)\sqrt\tau + \sqrt{G(\mu\tau)^2\tau - 2G(\mu\tau)\tau + 2}
\)
are \emph{simultaneously} verified for some $\mu>1$, with the function $G$ defined in \eqref{def:G}.
Then $\E[w(t)^2]$ decays exponentially to zero as $t\to\infty$,
where $w=w(t)$ is the solution of \eqref{eq:main} 
subject to the constant initial datum $w_0\neq 0$.
\end{theorem}

\begin{figure}
\centerline{
\includegraphics[width=0.6\columnwidth]{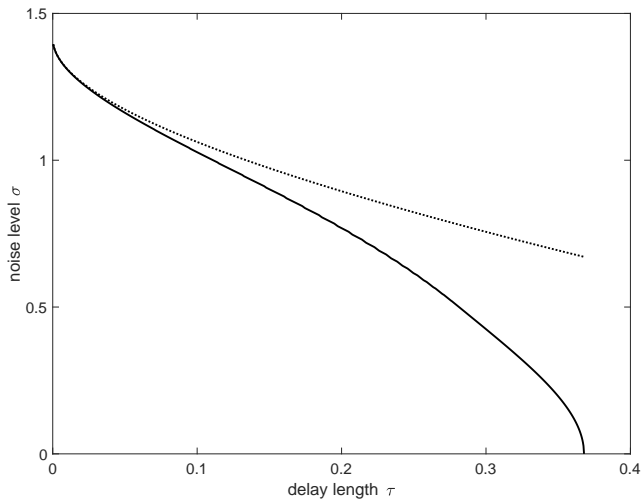}}
\caption{Numerical realization of the sufficient condition for exponential decay of $\E[w^2]$
as formulated in Theorem \ref{thm:exponential}, solid line. For comparison,
sufficient condition for asymptotic decay given by Theorem \ref{thm:asymptotic}, dotted line.
}
\label{fig:2}
\end{figure}

Obviously, the condition posed by Theorem \ref{thm:exponential} is not explicit, since it involves
a search for $\mu>1$ such that both \eqref{ass:sigma1} and \eqref{ass:sigma2} are satisfied.
Finding the maximal admissible $\sigma=\sigma(\tau)$ for a given $\tau\geq 0$ in fact means
\(  \label{maxmin}
   \sigma(\tau) := \max_{\mu>1} \min\left\{ f_1(\tau,\mu), f_2(\tau,\mu) \right\},
\)
where $f_1=f_1(\tau,\mu)$ and, resp., $f_2=f_2(\tau,\mu)$ denote
the right-hand sides of \eqref{ass:sigma1} and, resp., \eqref{ass:sigma2}.
It does not seem feasible to find an explicit analytical formula for $\sigma(\tau)$
in \eqref{maxmin}, however, the problem is quite easily approachable numerically.
First, let us observe that \eqref{ass:sigma1} is only satisfiable if $\mu>e^{\mu\tau}$,
which requires $\tau<e^{-1}$.
Consequently, for each $\tau\in(0,e^{-1})$ we only need to search values of $\mu$
such that $\mu>e^{\mu\tau}$, which represents a bounded interval.
The situation is also simplified by the fact that,
as revealed by a simple analysis, $f_2(\tau,\mu$)
is a decreasing function of $\mu$ for any fixed $\tau<e^{-1}$. 
The result of numerical realization of \eqref{maxmin} is plotted in Fig. \ref{fig:2},
where also the condition for asymptotic decay \eqref{cond:asympt2} is indicated for comparison.
Finally, let us note that for $\tau=0$, the conditions \eqref{ass:sigma1}--\eqref{ass:sigma2}
collapse to $\sigma^2<2$, which is the sharp condition for asymptotic decay
in mean square of the (nondelay) geometric Brownian motion.
Also the condition $\tau<e^{-1}$ is sharp, since all nontrivial solutions
of \eqref{eq:u} oscillate if $\tau>e^{-1}$.

\section{Asymptotic decay: Proof of Theorem \ref{thm:asymptotic}}\label{sec:asymptotic}
For $p$, $q>0$ and $t>0$ we define the functional
\(  \label{Layp}
   \L(t) := |w(t)|^2 + q \int_{t-\tau}^t |w(s)|^2 \d s + p \int_{t-\tau}^t \int_\theta^t |w(s-\tau)|^2 \d s\d\theta,
\)
where $w=w(t)$ is the solution of \eqref{eq:main}--\eqref{IC:main};
we refer to \cite{Fridman, Sun-Chen} for an overview of the theory of Lyapunov functionals
for systems with delay.

\begin{lemma}\label{lem:asDecay}
Let $\sigma^2<2$ and
\(  \label{cond:decay}
   \tau < 1 - \frac{\sigma}{2} \sqrt{4-\sigma^2}.
\)
Then there exist $p$, $q>0$ and $\kappa>0$ such that
\(  \label{decay}
   \tot{}{t} \E[\L(t)] \leq - \kappa\E[w(t-\tau)^2] \qquad\mbox{for } t>\tau.
\)
\end{lemma}

\proof
We apply the It\^{o} formula to calculate $\d |w(t)|^2$. 
Note that the It\^{o} formula holds in its usual form 
also for delay stochastic processes, see page 32 in \cite{Gillouzic}
or \cite{KM92, KM99, EN73, Mao2}, and with \eqref{eq:main} it gives
\(   \label{dw2}
    \d |w(t)|^2 = 2 \left(-\tw w + \frac{\sigma^2}{2} \tw^2 \right) \d t + 2\sigma\tw w \d B_t.
\)
Consequently,
\(  \label{dLayp}
   \d\L(t) &=& 2 \left(-\tw w + \frac{\sigma^2}{2} \tw^2 \right) \d t + 2\sigma\tw w \d B_t \\
      &&\qquad +\; q(w^2 - \tw^2)\d t + p\left( - \int_{t-\tau}^t |w(s-\tau)|^2 \d s + \tau \tw^2 \right)\d t.
    \nonumber
\)
For any $\delta>0$ we have
\[
   -2 \tw w = -2(\tw-w)w - 2w^2 \leq \delta |w-\tw|^2 + (\delta^{-1}-2) w^2.
\]
Restricting to $t>\tau$, we have for any $\eps>0$, 
\[
   |w-\tw|^2 = \left| \int_{t-\tau}^t \d w(s) \right|^2 &\leq& \left( \left| \int_{t-\tau}^t w(s-\tau) \d s \right| + \sigma \left| \int_{t-\tau}^t w(s-\tau) \d B_s \right| \right)^2 \\
      &\leq& (1+\eps) \left( \int_{t-\tau}^t w(s-\tau) \d s \right)^2 + (1+\eps^{-1}) \sigma^2 \left( \int_{t-\tau}^t w(s-\tau) \d B_s \right)^2.
\]
We take the expectation and use the Jensen inequality and Fubini theorem for the term
\[
   \E \left[ \left( \int_{t-\tau}^t w(s-\tau) \d s \right)^2 \right] \leq \tau \E \left[ \int_{t-\tau}^t |w(s-\tau)|^2 \d s \right]
      =  \tau \int_{t-\tau}^t \E\left[ |w(s-\tau)|^2 \right] \d s,
\]
and the isometry of the It\^{o} integral \cite[Theorem 5.8(iii)]{Mao}, 
for the term
\[
   \E \left[ \left( \int_{t-\tau}^t w(s-\tau) \d B_s \right)^2 \right] = \int_{t-\tau}^t \E\left[ |w(s-\tau)|^2 \right] \d s.
\]
Therefore, we arrive at
\[
   \E \left[ |w-\tw|^2 \right] \leq \left( (1+\eps)\tau + (1+\eps^{-1}) \sigma^2 \right) \int_{t-\tau}^t \E\left[ |w(s-\tau)|^2 \right] \d s.
\]
Minimization of the right-hand side in $\eps>0$ leads to $\eps:=\sigma/\sqrt{\tau}$, and thus
\[
   \E \left[ |w-\tw|^2 \right] \leq \left( \sqrt{\tau} + \sigma \right)^2 \int_{t-\tau}^t \E\left[ |w(s-\tau)|^2 \right] \d s.
\]
Consequently, taking the expectation in \eqref{dLayp}, we obtain
\[
   \tot{}{t} \E[\L(t)] \leq \left( \delta^{-1}+q-2\right) \E[w^2] + \left(\sigma^2 + p\tau - q\right) \E[\tw^2]
      + \left( \delta(\sqrt{\tau} + \sigma)^2 - p \right) \int_{t-\tau}^t \E\left[ |w(s-\tau)|^2 \right] \d s.
\]
With the choice
\[
    p:=\delta(\sqrt{\tau} + \sigma)^2, \qquad q:=2-\delta^{-1}
\]
we arrive at
\[
   \tot{}{t} \E[\L(t)] \leq \Bigl(\sigma^2 + \delta \tau(\sqrt{\tau} + \sigma)^2 + \delta^{-1} - 2\Bigr) \E[\tw^2].
\]
Minimization of the right-hand side with respect to $\delta>0$ gives $\delta:=\left(\tau(\sqrt{\tau} + \sigma)^2\right)^{-1/2}$ and
\[
   \tot{}{t} \E[\L(t)] \leq -\kappa \E[\tw^2]
\]
with
\(   \label{kappa}
   -\kappa := \left(\sigma+\sqrt\tau\right)^2 + \tau - 2.
\)
Finally, a simple calculation reveals that if $\sigma^2 <2$, then $\kappa > 0$
if and only if \eqref{cond:decay} is satisfied.
It is easily checked that then $p$ and $q$ are both positive.
\endproof

\textbf{Proof of Theorem \ref{thm:asymptotic}.}
Obviously, with the bounded initial datum \eqref{IC:main}, we have $\E[\L(\tau)] < +\infty$
due to \eqref{dLayp}.
An integration of \eqref{decay} in time gives, for $t>\tau$,
\(
   \E[w(t)^2] &\leq& \E[\L(t)] = \E[\L(\tau)] + \int_{\tau}^t \tot{}{s} \E[\L(s)] \d s  \nonumber \\
      &\leq& \E[\L(\tau)] - \kappa \int_{\tau}^t \E[w(s-\tau)^2] \d s,
      \label{intConv}
\)
with $\kappa>0$ given by \eqref{kappa}.
Consequently, $\E[w(t)^2]$ is uniformly bounded by $\E[\L(\tau)]$ for $t>\tau$.
Taking the expectation in \eqref{dw2} and using the Cauchy-Schwartz inequality, we have
\[
   \tot{}{t} \E[w(t)^2] &=& -2 \E[w\tw] + \sigma^2 \E[\tw^2] \\
      &\leq& \E[w^2] + (1+\sigma^2) \E[\tw^2],
\]
which gives uniform boundedness of $\tot{}{t} \E[w(t)^2]$ for $t>2\tau$.
Moreover, we note that due to \eqref{intConv}
the integral $\int_{\tau}^t \E[w(s-\tau)^2] \d s$ is convergent as $t\to\infty$.
Barbalat's lemma \cite{Barb} then implies that $\lim_{t\to\infty} \E[w(t)^2] = 0$
and concludes the proof of Theorem \ref{thm:asymptotic}.

\section{Exponential decay: Proof of Theorem \ref{thm:exponential}}\label{sec:exponential}
In this section we assume that $w=w(t)$ is a solution of \eqref{eq:main} 
subject to the deterministic constant initial datum $w_0\neq 0$, and we introduce the notation
\(
   y(t) &:=& \E[w^2(t)/2] \qquad\mbox{for } t\geq 0,\\
         &:=& w_0^2/2 \qquad\mbox{for } t< 0.     \nonumber
\)

\begin{lemma}\label{lem:1}
Let $\sigma^2\leq 2$.
If for some $\mu>1$ the condition
\(  \label{cond:mu1}
    2e^{\mu\tau} + \sigma^2 e^{2\mu\tau} \leq 2\mu
\)
is satisfied, then for all $t\in\R$ and $s>0$,
\(  \label{claim:lem:1}
   e^{-2\mu s} y(t) < y(t-s) < e^{2\mu s} y(t).
 \)
\end{lemma}

\proof
An application of the It\^{o} formula gives
\[
   \frac{\d{w^2}}{2} = \left(-\tw w + \frac{\sigma^2}{2} \tw^2 \right) \d t + \sigma\tw w \d B_t,
\]
and taking expectation, we have for $t>0$,
\(  \label{toty}
   \dot y = \E\left[ -\tw w + \frac{\sigma^2}{2} \tw^2 \right].
\)
With the constant initial datum $w_0\in\R$, we obtain for $t=0$,
\[
   \dot y(0+) = \left( - 1 + \frac{\sigma^2}{2} \right) w_0^2,
\]
where $\dot y(0+)$ denotes the right-hand side derivative of $y$ at $t=0$.
Consequently, since by assumption $\mu>1$,
\[
    \left| \frac{\dot y(0+)}{y(0)} \right| = 2-\sigma^2 < 2\mu.
\]
Due to nonzero constant initial datum
and the continuity of $y(t)$ and $\dot y(t)$ for $t>0$, 
there exists $T>0$ such that
\(   \label{uptoT}
    \left| \frac{\dot y(t)}{y(t)} \right| < 2\mu \qquad\mbox{for } t<T.
\)
We claim that \eqref{uptoT} holds for all $t\in\R$, i.e., that $T=+\infty$.

For contradiction, assume that $T<+\infty$, then again by continuity we have
\( \label{uptoT2}
   |\dot{y}(T)| = 2\mu y(T).
\)
Integrating \eqref{uptoT} on the time interval $(T-s,T)$ with $s>0$ yields
\(  \label{est_temp}
      y(T-s) < e^{2\mu s} y(T).
\)
With \eqref{toty} and the Cauchy-Schwartz inequality with $\eps>0$ we have for $t>0$,
\[
   |\dot y| \leq \E\left[ |\tw| |w|\right] + \sigma^2 \ty \leq 2\sqrt{\ty}\sqrt{y} + \sigma^2 \ty
      \leq \eps y + \left(\eps^{-1} + \sigma^2\right) \ty.
\]
Using \eqref{est_temp} with $s:=\tau$ gives $y(T-\tau) < e^{2\mu\tau} y(T)$, so that
\[
   |\dot y(T)| < \Bigl(\eps + \left(\eps^{-1} + \sigma^2\right) e^{2\mu\tau}\Bigr) y(T),
\]
and minimization of the right-hand side with respect to $\eps>0$ gives
\[
   |\dot y(T)| < \left(2e^{\mu\tau} + \sigma^2 e^{2\mu\tau}\right) y(T).
\]
Finally, assumption \eqref{cond:mu1} gives
\[
   |\dot y(T)| < 2\mu y(T),
\]
which is a contradiction to \eqref{uptoT2}.
Consequently, \eqref{uptoT} holds for all $t\in\R$,
and an integration on the interval $(t-s,t)$,
taking into account the constant initial datum, implies \eqref{claim:lem:1}.
\endproof

\begin{lemma} \label{lem:2}
Let the condition \eqref{cond:mu1} of Lemma \ref{lem:1} be satisfied for some $\mu>1$.
Then we have, along the solutions of \eqref{eq:main},
\(   \label{claim:lem:decay}
   \dot y < 2 \left[ \left( \sqrt{\tau} + \sigma \right) \left( \frac{e^{2\mu\tau}-1}{2\mu} \right)^{1/2} + \frac{\sigma^2}{2} - 1 \right] \ty,
\)
for $t>0$.
\end{lemma}

\proof
Referring to \eqref{toty} we have for $t>0$,
\(  \label{toty2}
   \dot y = \E\left[ -\tw w \right] + \frac{\sigma^2}{2} \E\left[ \tw^2 \right],
\)
and, with the Cauchy-Schwartz inequality,
\[
   \E\left[ -\tw w \right] = \E\left[(\tw-w)\tw \right] - \E\left[\tw^2\right]
      \leq \Bigl( \E\left[ |w-\tw|^2 \right] \Bigr)^{1/2} \Bigl( \E \left[\tw^2\right] \Bigr)^{1/2} - \E \left[\tw^2\right].
\]
If $t>\tau$, we have for any $\eps>0$, 
\[
   |w-\tw|^2 = \left| \int_{t-\tau}^t \d w(s) \right|^2 &=& \left( \int_{t-\tau}^t w(s-\tau) \d s  + \sigma  \int_{t-\tau}^t w(s-\tau) \d B_s \right)^2 \\
      &\leq& (1+\eps) \left( \int_{t-\tau}^t w(s-\tau) \d s \right)^2 + (1+\eps^{-1}) \sigma^2 \left( \int_{t-\tau}^t w(s-\tau) \d B_s \right)^2.
\]
As in the proof of Lemma \ref{lem:asDecay}, we take the expectation and use the Jensen inequality,
Fubini theorem and isometry of the It\^{o} integral to obtain
\(   \label{ineq}
   \E \left[ |w-\tw|^2 \right] \leq \left( (1+\eps)\tau + (1+\eps^{-1}) \sigma^2 \right) \int_{t-\tau}^t \E\left[ |w(s-\tau)|^2 \right] \d s.
\)
If $0 < t \leq\tau$, we have, due to the constant initial condition,
\[
   |w-\tw|^2 = \left| \int_{t-\tau}^t \d w(s) \right|^2 = \left( \int_{0}^t w(s-\tau) \d s  + \sigma  \int_{0}^t w(s-\tau) \d B_s \right)^2,
\]
and a trivial modification of the above estimates gives \eqref{ineq} again.
Consequently,  \eqref{ineq} holds for all $t>0$,
with the constant initial datum $w_0$ being extended to the interval $[-2\tau,0]$.
Minimization of the right-hand side in $\eps>0$ leads to $\eps:=\sigma/\sqrt{\tau}$, and thus
\[
   \E \left[ |w-\tw|^2 \right] \leq \left( \sqrt{\tau} + \sigma \right)^2 \int_{t-\tau}^t \E\left[ |w(s-\tau)|^2 \right] \d s.
\]
An application of Lemma \ref{lem:1} gives
\[
   \int_{t-\tau}^t \E\left[ |w(s-\tau)|^2 \right] \d s = 2 \int_{t-\tau}^t y(s-\tau) \d s < 2 y(t-\tau) \int_0^\tau e^{2\mu s} \d s = \frac{e^{2\mu\tau}-1}{\mu} \ty.
\]
Consequently, we have
\[
   \E\left[ -\tw w \right] < 2 \left[ \left( \sqrt{\tau} + \sigma \right) \left( \frac{e^{2\mu\tau}-1}{2\mu} \right)^{1/2} - 1 \right] \ty,
\]
and inserting this into \eqref{toty2} immediately gives \eqref{claim:lem:decay}.
\endproof

\textbf{Proof of Theorem \ref{thm:exponential}.}
Lemmata \ref{lem:1} and \ref{lem:2} assert that
$y(t) = \E[w(t)^2/2]$ is monotonically decaying if condition
\eqref{cond:mu1} is satisfied for some $\mu>1$ and if
\(   \label{cond:mu2}
    \left( \sqrt{\tau} + \sigma \right) \left( \frac{e^{2\mu\tau}-1}{2\mu} \right)^{1/2} + \frac{\sigma^2}{2} < 1.
\)
A simple calculation reveals that \eqref{cond:mu1} is equivalent to \eqref{ass:sigma1}, 
while \eqref{cond:mu2} is equivalent to \eqref{ass:sigma2}.
Then, a combination of \eqref{claim:lem:1} with $s:=\tau$ and \eqref{claim:lem:decay} yields
\[
   \dot y(t) < \left[ \left( \sqrt{\tau} + \sigma \right) \left( \frac{e^{2\mu\tau}-1}{2\mu} \right)^{1/2} + \frac{\sigma^2}{2} - 1 \right] e^{-2\mu\tau} y(t),
\]
for $t>0$, which implies exponential decay of $y=y(t)$ in time.

{\small
{\em Authors' addresses}:
{\em Jan Haskovec}, King Abdullah University of Science and Technology, 23955 Thuwal, KSA.
e-mail: \texttt{jan.haskovec@\allowbreak kaust.edu.sa}.

}

\end{document}